\documentclass[11pt,a4paper]{article}
\usepackage{epsfig}
\usepackage{amssymb}
\usepackage{epic}
\usepackage{eepic}
\usepackage[latin1]{inputenc}
\usepackage[T1]{fontenc}
\usepackage{psfrag}
\def\P{{\mathbb P}}

\def\Cox{\hfill \Box}
\def\GM{{\rm GM}}
\def\supp{{\rm supp}}
\newtheorem{theorem}{Theorem}[section]
\newtheorem{lemma}[theorem]{Lemma}
\newtheorem{corollary}[theorem]{Corollary}
\newtheorem{proposition}[theorem]{Proposition}
\newtheorem{definition}[theorem]{Definition}
\newtheorem{example}[theorem]{Example}
\newtheorem{remark}[theorem]{Remark}

\author{Erik I. Broman\footnote{Dept.\ of Mathematics,
Chalmers University of Technology, S-412 96 G\"oteborg, Sweden, 
\texttt{http://www.math.chalmers.se/\~{ }broman/}, research 
supported by the Swedish Research Council.}, 
Olle H\"aggstr\"om\footnote{Dept.\ of Mathematics,
Chalmers University of Technology, S-412 96 G\"oteborg, Sweden, 
\texttt{http://www.math.chalmers.se/\~{ }olleh/}, research 
supported by the Swedish Research Council.}, and 
Jeffrey E. Steif\footnote{Dept.\ of Mathematics,
Chalmers University of Technology, S-412 96 G\"oteborg, Sweden, 
\texttt{http://www.math.chalmers.se/\~{ }steif/}, research 
supported by the Swedish Research Council
and the G\"{o}ran Gustafsson Foundation (KVA).}}
\date{\today}

\title{Refinements of stochastic domination}

\begin{document}
  \maketitle

  \begin{abstract}
In a recent paper by two of the authors, 
the concepts of upwards and downwards 
$\epsilon$-movability were introduced, mainly as a 
technical tool for studying dynamical percolation of 
interacting particle systems. 
In this paper, we further explore these concepts which can be seen as 
refinements or quantifications of stochastic domination, and we
relate them to previously studied concepts such as uniform 
insertion tolerance and extractability.



%  \medskip\noindent
%  {\bf Short title:} Determinants and two-block-factors

  \end{abstract}

  \section{Introduction}\label{intro}

  In \cite{BS}, Broman and Steif introduced certain refinements
of stochastic domination, which we call upwards and 
downwards  
$\epsilon-$movability; see Definition \ref{def_eps-stab} below. 
These concepts were introduced mainly as a technical tool in the analysis 
of dynamical percolation for interacting particle systems, but 
they turn out to be interesting in their own right.
In the present paper, we 
explore these concepts further and relate them to various previously
studied concepts.

Let $S$ be a countable set. For $p \in [0,1]$, 
let every $s\in S$, independently of all other elements in $S$,
take value $1$ with probability $p$ and take value $0$ 
with probability $1-p$.  We let $\pi_{p}$ 
denote the corresponding product measure on $\{0,1\}^{S}$. 
When talking about product measures on $\{0,1\}^S$, 
we will always mean these uniform ones (with
the same $p$ for every $s \in S$). 

Let $\mu$ be an arbitrary 
probability measure on $\{0,1\}^{S}$. For $\epsilon \in (0,1)$,
we will let $\mu^{(+,\epsilon)}$ denote the distribution of the process
obtained by first choosing an element of $\{0,1\}^{S}$ according to
$\mu$ and then independently changing each 0 to a 1 with probability
$\epsilon$. Similarly, we will let $\mu^{(-,\epsilon)}$ denote the 
distribution of the process obtained by first choosing an element of 
$\{0,1\}^{S}$ according to $\mu$ and then independently changing each 1 to 
a 0 with probability $\epsilon$. Finally, for $\delta \in (0,1)$,
we let $\mu^{(-,\epsilon,+,\delta)}$
denote the distribution of the process obtained by first choosing an element 
of $\{0,1\}^{S}$ according to $\mu$ and then independently changing each 0 to 
a 1 with probability $\delta$ and each 1 to a 0 with probability $\epsilon$.

It turns out that for any $\epsilon\in [0,1)$,
$\mu_1^{(+,\epsilon)}=\mu_2^{(+,\epsilon)}$ implies 
that $\mu_1=\mu_2$. To
see this, it suffices to check that $\mu_1(A)=\mu_2(A)$ for events $A$ of 
the type ``all $s \in S'$ take value $0$'' where $S' \subseteq S$
(this is easy),
and then use inclusion-exclusion. Similarly, of course,  
$\mu_1^{(-,\epsilon)}=\mu_2^{(-,\epsilon)}$ 
implies $\mu_1=\mu_2$.

  For $\sigma,\sigma^{\prime}\in\{0,1\}^{S}$ we write 
$\sigma \preceq \sigma^{\prime}$
  if $\sigma(s)\leq \sigma^{\prime}(s)$ for every $s\in S.$ A
  function $f:\{0,1\}^{S} \rightarrow {\mathbb R}$ 
is increasing if $f(\sigma)\leq f(\sigma^{\prime})$ whenever 
$\sigma \preceq \sigma^{\prime}.$ 
  For two probability measures $\mu,\mu^{\prime}$ on $\{0,1\}^{S}$, we
say that $\mu$ {\bf is stochastically dominated by} $\mu'$, and
  write $\mu \preceq \mu^{\prime}$, if for every continuous
  increasing function $f$ we have that $\mu(f)\leq \mu^{\prime}(f)$.
($\mu(f)$ is shorthand for $\int f d\mu$.) By Strassens theorem
(see \cite[p.\ 72]{IPS}), 
this is equivalent to the existence of random variables
$X,X^{\prime}\in\{0,1\}^{S}$ such that 
$X$ has distribution $\mu$, $X^{\prime}$ 
has distribution $\mu^{\prime}$, and $X \preceq X^{\prime}$ a.s.
This is also equivalent to $\mu(A)\leq \mu^{\prime}(A)$ for all
up-sets $A$ where an up-set is a set whose indicaor function 
$I_A$ is increasing. From now on ``$\sim$'' will mean ``has distribution''.

  \begin{definition} \label{def_eps-stab} 
  Let $(\mu_{1},\mu_{2})$ be a pair of 
probability measures on $\{0,1\}^{S},$ where $S$ is a countable set. 
  Assume that $\mu_{1} \preceq \mu_{2}$.
  If, given $\epsilon>0$, we have
  \[
    \mu_{1} \preceq  \mu_{2}^{(-,\epsilon)},
  \]
  then we say that the pair $(\mu_{1},\mu_{2})$ is
  {\bf downwards $\epsilon$-movable}. 
$(\mu_{1},\mu_{2})$ is said to be {\bf downwards movable}
if it is downwards 
$\epsilon$-movable for some $\epsilon >0$.
Analogously, if, given $\epsilon>0$, we have 
  \[
    \mu_{1}^{(+,\epsilon)} \preceq  \mu_{2},
  \] 
  then we say that the pair ($\mu_{1},\mu_{2}$) is 
{\bf upwards $\epsilon$-movable},
and we say that ($\mu_{1},\mu_{2}$)
is {\bf upwards movable} if the pair is upwards $\epsilon$-movable 
for some $\epsilon >0$.
  \end{definition}

Note that if we restrict to the case where both $\mu_1$ and $\mu_2$ 
are product measures, then these concepts become trivial.

In \cite{BS} a considerable amount of effort was spent on trying to show 
downwards movability when the pair
considered was two stationary 
distributions, corresponding to two different parameter values,
for some specific interacting particle system. 
In particular, the so called contact process (see Liggett \cite{SIS}
for definitions and a survey) was investigated.
Considering $(\mu_1,\mu_2)$, 
where $\mu_i$ is the upper invariant measure for the contact process
with infection rate $\lambda_i$,
it was shown in \cite{BS}
that if $\lambda_1<\lambda_2$, then the pair is downwards movable. 

Another result from \cite{BS} 
is that if $\mu_1 \preceq \mu_2$, 
$\mu_2$ satisfies the FKG lattice condition (see \cite[p.\ 78]{IPS})
and 
  \[
    \displaystyle \inf_{\tilde{S} \subset S \atop |\tilde{S}|<\infty} 
    \inf_{s\in \tilde{S}\atop \xi \in \{0,1\}^{\tilde{S}\setminus s}}
    [\mu_{2}(\sigma(s)=1|\sigma(\tilde{S}\setminus s)\equiv \xi)-
     \mu_{1}(\sigma(s)=1|\sigma(\tilde{S}\setminus s)\equiv \xi)] 
>0
  \]
then $(\mu_1,\mu_2)$ is downwards movable. 
This  however is not sufficient to 
get the result for the contact process mentioned above since
by \cite{L1}, the upper invariant measure for the contact process on
${\mathbb Z}$ does not satisfy the FKG lattice condition when 
$\lambda < 2$.

In the present paper, we will concentrate on the case 
where $\mu_{1}$ is a product measure but $\mu_{2}$ is not.
We now proceed with some further explanations and definitions needed  
to state our main results, Theorems \ref{thm1} and \ref{thm:finite}
below.  
In Sections \ref{examples}--\ref{secFKG}, 
we will establish a number of examples and auxiliary results, 
while Section \ref{types} will tie things together giving proofs
of Theorems \ref{thm1} and \ref{thm:finite}.

For a probability measure $\mu$ on $\{0,1\}^{S}$,  
  define $p_{\sup,\mu}$ by 
  \[
    p_{\sup,\mu}:=\sup \{p \in [0,1]: \pi_{p} \preceq \mu \}.
  \]
  Since the relation $\preceq$ is preserved under weak limits we see that
  \[
    \pi_{p_{\sup,\mu}} \preceq\mu
  \]
  and so the supremum is achieved. Therefore we also denote this by
 $p_{\max,\mu}$.

If $p_{\max,\mu}=0$, then trivially $(\pi_{p_{\max,\mu}},\mu)$
is downwards movable but not upwards movable. Assume next
that $\mu$ is a probability measure with $p_{\max,\mu}>0$.
If $p \in [0,p_{\max,\mu})$, 
then the pair $(\pi_{p},\mu)$ is trivially upwards
movable. It is also easy to see that it is downwards movable 
by arguing as follows.
By Strassen's theorem, we may choose
$X \sim \mu$ and $Y \sim \pi_{p_{\max,\mu}}$ such that $X \geq Y$ a.s.
Then choose $\epsilon>0$ such that
$(1-\epsilon)p_{\max,\mu}>p$, and let $Z \sim \pi_{1-\epsilon}$ 
be independent of both $X$ and $Y$. 
We obtain $\min(X,Z) \geq \min(Y,Z)$ a.s., and
since $\min(Y,Z) \sim \pi_{p_{\max,\mu}(1-\epsilon)}$ we conclude that
  \[
    \pi_{p} \preceq \mu^{(-,\epsilon)},
  \]
as desired. 

The final case we are left with (when one of the measures is a uniform
product measure) is $(\pi_{p_{\max,\mu}},\mu)$ with $p_{\max,\mu}>0$.
This pair is by definition not upwards movable, 
but we believe it is an interesting
question to ask if it is downwards movable and 
this question motivates the following definition.
  \begin{definition}
  We say that $\mu$ is {\bf nonrigid} if the pair
$(\pi_{p_{\max,\mu}},\mu)$ is downwards movable and
otherwise we will say that $\mu$ is {\bf rigid}. 
  \end{definition}

All uniform product measures other than $\delta_0$ are trivially rigid
while all $\mu$ such that $p_{\max,\mu}=0$ are trivially nonrigid.
Heuristically, it is natural to expect that as long as
$p_{\max,\mu}>0$, then typically $\mu$ should be rigid.
This issue turns out to be quite intricate, however; see
Proposition \ref{prop:rigid_when_finite} and Theorem 
\ref{extype2_2} below. 

Other well known concepts which have arisen in a number of
problems and which we feel belong to this same circle of
ideas are those of finite energy (Newman and Schulman \cite{NS})
and insertion and deletion tolerance (Lyons and Schramm \cite{LSch}).

  \begin{definition}
  We say that $\mu$ is {\bf $\epsilon$-insertion tolerant} if 
for any $s\in S$, we have that
  \begin{equation} \label{eqnUIT}
    \mu(\sigma(s)=1|\sigma(S \setminus s))\geq 
\epsilon \textrm{ a.s. }
  \end{equation}
We say that $\mu$ is {\bf uniformly insertion tolerant} if it is
 $\epsilon$-insertion tolerant for some $\epsilon>0$. The analogous notions of
{\bf $\epsilon$-deletion tolerant} and {\bf uniformly deletion tolerant} are
defined similarly (the ``$1$'' is replaced by ``$0$'').
Finally, we say $\mu$ has {\bf finite $\epsilon$-energy} if it is both
$\epsilon$-insertion tolerant and $\epsilon$-deletion tolerant, and that it has
{\bf uniform finite energy} 
if it has finite $\epsilon$-energy for some $\epsilon>0$.
  \end{definition}
Also closely related are the following notions of extractability; 
we discuss some background on
this concept at the end of the introduction.
\begin{definition}
We call
$\mu$ {\bf $\epsilon$-upwards extractable} if
there exists a probability measure $\nu$ 
such that $\mu=\nu^{(+,\epsilon)}$. 
We call $\mu$ {\bf uniformly upwards extractable} if
it is $\epsilon$-upwards extractable for some $\epsilon>0$. 
The notions of {\bf $\epsilon$-downwards extractable} and
{\bf uniformly downwards extractable} are defined analogously
(the ``$+$'' is replaced by ``$-$''). Finally,
$\mu$ is called {\bf $\epsilon$-extractable} if
there exists a probability measure $\nu$ 
such that $\mu=\nu^{(-,\epsilon,+,\epsilon)}$,
and it is called {\bf uniformly extractable} if it is
$\epsilon$-extractable for some $\epsilon>0$.
  \end{definition}
We are now equipped with all the definitions needed to state our main 
theorem. We refer to Figure 1 for a comprehensive diagram over the 
implications and non-implications that the theorem asserts.

  \begin{theorem} \label{thm1}
  Let $S$ be a countable set and consider the following 
properties of a probability measure $\mu$ on 
  $\{0,1\}^S$:
  \begin{description}
  \item{\rm (I) } $\mu$ is uniformly upwards extractable.
  \item{\rm (II) } $\mu$ is uniformly insertion tolerant.
  \item{\rm (III) } $\mu$ is rigid.
  \item{\rm (IV) } There exists a $p>0$ such that $\pi_p \preceq \mu.$
  \end{description}
We then have that {\rm (I)} $\Rightarrow$ {\rm (II)} $\Rightarrow$ {\rm (IV)}
and that {\rm (I)} $\Rightarrow$ {\rm (III)} $\Rightarrow$ {\rm (IV)} while 
none of the four corresponding reverse implications hold. Also,
{\rm (III)} does not imply {\rm (II)}. 
Moreover, with $S={\mathbb Z}$, there exist translation invariant 
examples for all of the asserted nonimplications. 
%  \[
%  {\rm (I)}
%   \left. \begin{array}{c} \Rightarrow   \end{array} \right. 
%  {\rm (II)}
%   \left. \begin{array}{c} \Rightarrow \\ \not \Leftarrow  \end{array} \right. 
%  {\rm (IV)}
%  \]
%  \[
%{\rm (I)}
%   \left. \begin{array}{c} \Rightarrow \\ \not \Leftarrow  \end{array} \right. 
%  {\rm (III)}
%   \left. \begin{array}{c} \Rightarrow \\ \not \Leftarrow  \end{array} \right. 
%  {\rm (IV)}
%  \]
%  \[
%{\rm (II)}
%   \left. \begin{array}{c} \not \Rightarrow \\ \not \Leftarrow  \end{array} \right. 
%{\rm (III)}
%  \]
%  \[
%{\rm (IV)}
%   \left. \begin{array}{c} \not \Rightarrow  \end{array} \right. 
%{\rm (II)}  \mbox{ or } {\rm (III)}
%  \]
%
  \end{theorem}
In addition, it turns out that {\rm (IV)} does not even imply 
``{\rm (II)} or {\rm (III)}''; see Remark \ref{rem:neither(II)or(III)}.
Note that we have not managed to work out whether or not (II) implies (III).

\begin{figure}[hbt]
\begin{center}
\psfrag{1}{I}
\psfrag{2}{\hspace{-1mm}II}
\psfrag{3}{\hspace{-1mm}III}
\psfrag{4}{\hspace{-1mm}IV}
\includegraphics[width=0.4\textwidth]{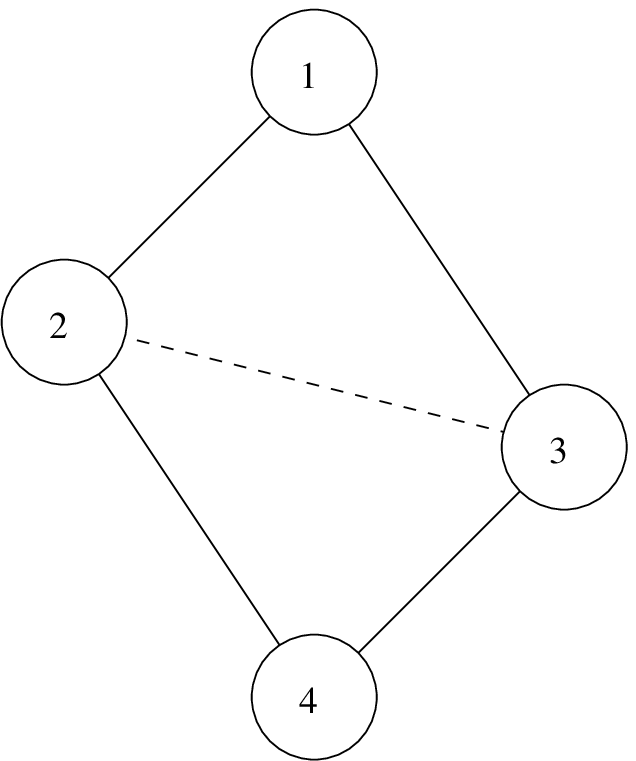}
\end{center}
\end{figure}

\begin{quote}
{\small Figure 1. 
Hasse diagram of the implications between properties (I), (II),
(III) and (IV) in Theorem \ref{thm1}: we have proved that
one property implies another iff there
is a downwards path in the diagram from the former to the latter. We do not
know whether the dashed line between (II) and (III) should be there or not,
i.e., whether or not uniform insertion tolerance implies rigidity. 
As will be seen in Theorem \ref{thmFKG}, the
desired implication (II) $\Rightarrow$ (III) holds under an additional
FKG-like assumption. If we restrict to finite $S$, then
some of the implications will turn into equivalences; see Theorem 
\ref{thm:finite}.}
\end{quote}

Some of the asserted implications are easy: (I) trivially implies (II).
The implication (III) $\Rightarrow $ (IV) is also trivial as we saw.
It is a direct application of Holley's inequality (see, e.g., 
\cite[Theorem 4.8]{GHM}) to see that
$\epsilon$-insertion tolerance implies that $\pi_\epsilon \preceq \mu$,
whence (II) implies (IV).
Thus, apart from the implication (I) $\Rightarrow$ (III) 
(which is in fact not so hard either), we see all the implications claimed
in the theorem. Therefore
our interest in Theorem \ref{thm1} is more in the counterexamples showing
the distinction between some of these properties rather than in the 
implications.

As mentioned above, we do not know in general whether (II) implies (III). 
However Theorem \ref{thmFKG} provides us 
with a partial answer, telling us that this is true under 
the extra assumption of $\mu$ being downwards FKG,
a property weaker than satisfying the FKG lattice condition and
defined as follows.
  \begin{definition} \label{downFKG}
  A measure $\mu$ on $\{0,1\}^S$ is downwards FKG if for any finite 
$S^{\prime}\subset S$ and any increasing subsets $A,B$
  \[
    \mu(A\cap B| \sigma(S^{\prime})\equiv 0)
\geq\mu(A | \sigma(S^{\prime})\equiv 0) \mu(B | \sigma(S^{\prime})\equiv 0).
  \]
  \end{definition}
The concept of downwards FKG was made explicit in \cite{LigS} where 
it was shown,
among other things, that for such translation invariant measures, stochastic
domination of a product measure has a very simple characterization. 
As mentioned before, the upper invariant measure for the contact process
in one dimension and with $\lambda<2$ is known to not satisfy the FKG lattice condition.
In addition this is believed to be true for any value of $\lambda$ and any dimension.
However, it was proven in \cite{BHK} that it is downwards FKG for any dimension
and for all values of $\lambda.$
This result was then
exploited in \cite{LigS} to show that the upper invariant measure for the 
contact process dominates product measures despite the fact that 
the measure is not uniformly insertion tolerant.

\begin{theorem} \label{thmFKG}
Let $\mu$ be a translation invariant downwards FKG
measure on $\{0,1\}^{{\mathbb Z}^{d}}.$ Then {\rm (II)} implies {\rm (III)}.
\end{theorem}
In Section \ref{secFKG} we prove an easy technical lemma that together 
with some results of \cite{LS} will give us the following theorem (see 
Section \ref{secFKG} for the definition of conditional negative associativity).

  \begin{proposition}\label{prop:neg}
Let $\mu$ be a translation invariant,
conditionally negatively associated measure on $\{0,1\}^{\mathbb Z}$.
Then {\rm (IV)} implies {\rm (III)}.
  \end{proposition}
If we restrict to finite $S$, then further implications between
the various properties are available. By the support of a 
measure $\mu$ on $\{0,1\}^S$, denoted by $\supp(\mu)$, we mean
$\{\xi \in \{0,1\}^S:\mu(\sigma(S)\equiv \xi)>0\}$.
\begin{theorem}  \label{thm:finite}
Let $S$ be finite, and consider properties {\rm (I)}--{\rm (IV)} of
probability measures on $\{0,1\}^S$. We then have
\begin{equation} \label{eqnSfin1}
{\rm (I)} \, \Leftrightarrow \, {\rm (II)} \Leftrightarrow \, \supp(\mu) \textrm{ is an up-set,}
\end{equation}
and
\begin{equation} \label{eqnSfin2}
{\rm (III)} \, \Leftrightarrow \, {\rm (IV)} \Leftrightarrow \, \mu(\sigma(S)\equiv 1)>0.
\end{equation}
Consequently, the properties in {\rm (\ref{eqnSfin1})} imply those in
{\rm (\ref{eqnSfin2})} 
but not vice versa. Note in particular that if 
we are in the full support case,
then {\rm (I)--(IV)} all hold.
\end{theorem} 
Although the term {\em extractability} is our own, the concept does 
have a history; in particular, there has been interest in finding
lower bounds
on $\epsilon$ for which $\epsilon$-extractability holds.
The question of uniform extractability has been
studied for the Ising model as well as other Markov random fields in 
\cite{BD,HB,N}. Earlier, in \cite{G1,G2,G3}, 
a similar question was studied for Markov chains and
autoregressive processes.
Of related interest is the result in \cite{HB} that for Markov random fields, 
uniform finite energy
implies uniform extractability.

\section{Basic examples} \label{examples}
  
Our first example in this section is a pair of measures which 
is downwards but not upwards movable. Note first that if 
$\nu^{(+,\epsilon)} \preceq  \mu$,   
then we must have 
\[
\nu \preceq  \mu
\]  
as well as
\[
\pi_{\epsilon} \preceq  \mu \, . 
\]  

\begin{example} \label{ex-epsstab1}
{\rm 
  Take $\nu=\frac{1}{2} \pi_{q}+\frac{1}{2}\delta_{0}$ and $\mu=\frac{1}{2} \pi_{p}+\frac{1}{2} \delta_{0}$ where 
  $q<p$ and where $\delta_{0}$ is the measure 
  which puts probability 1 on the configuration of all zeros. Trivially
  \[
    \nu \preceq  \mu.
  \]
If $S$ is infinite, then obviously
 $\mu$ cannot dominate a product measure with positive density.
  Therefore there does not exist any $\epsilon>0$ such that
  \[
    \nu^{(+,\epsilon)} \preceq  \mu.
  \]      
  However,
  \[
    \mu^{(-,\epsilon)}=
{\textstyle 
\frac{1}{2} \pi^{(-,\epsilon)}_{p}+\frac{1}{2} \delta^{(-,\epsilon)}_{0}
      =\frac{1}{2} \pi_{p(1-\epsilon)}+\frac{1}{2} \delta_{0}} \, ,
  \]
  so if we take $\epsilon>0$ such that $p(1-\epsilon)>q$, we get that
  \[
    \nu \preceq  \mu^{(-,\epsilon)}.
  \]
  Hence $(\nu,\mu)$ is downwards but not upwards movable.
} $\Cox$
\end{example}
  Before presenting the next three examples, we recall a family of
stationary processes known as {\bf determinantal} processes,
introduced in Lyons and Steif \cite{LS}. These are
  probability measures ${\mathbf P}^{f}$ on the Borel sets
  of $\{0,1\}^{{\mathbb Z}}$ where $f:[0,1] \rightarrow [0,1]$
  is a Lebesgue-measurable function (see \cite{LS}). For such
  an $f,$ define 
   \begin{eqnarray} \label{defPf} 
   \lefteqn{{\bf P}^{f}[\sigma(e_{1})=\cdots=\sigma(e_{k})=1]}\\
   & & :={\bf P}^{f}[\{\sigma \in \{0,1\}^{\mathbb Z}: \sigma(e_{1})=\cdots=\sigma(e_{k})=1\}]\nonumber \\
   & & :=\textrm{det}[\hat{f}(e_{j}-e_{i})]_{1\leq,i,j\leq k}, \nonumber 
   \end{eqnarray}
  where $e_{1},\ldots,e_{k}$ are distinct elements in ${\mathbb Z}$ 
  and $k\geq 1.$ Here $\hat{f}$ denotes the
  Fourier coefficients of $f,$ defined by
  \[
    \hat{f}(k):=\int_{0}^{1}f(x)e^{-i2\pi kx} dx.
  \]
  In \cite{LS} it is proven that ${\bf P}^{f}$ is
  indeed a probability measure. 
(The fact that a probability measure is determined by the values it gives
to cylinder sets of this type follows immediately from 
inclusion-exclusion.) In fact they showed this for the
  more general case of $f:{\mathbb T}^{d} \rightarrow [0,1]$ where
  ${\mathbb T}^{d}:={\mathbb R}^{d}/{\mathbb Z}^{d}$; in this case
  the resulting process is indexed by $ {\mathbb Z}^{d}$.
  This result rests very strongly on the
  results in \cite{Lyons2002}.  We will also need the following definition,
where GM stands for geometric mean:
  \[
    \GM(f):=\exp{\int_{0}^{1} \log f(x) dx}.
  \]

  \begin{example} \label{ex-epsstab2}
{\rm
  Let $f$ be a function from $[0,1]$ to itself. 
  By \cite[Theorem 5.3]{LS},  
$\pi_{p} \preceq {\mathbf P}^{f}$ iff $p\leq \GM(f)$. 
  It is easy to see from (\ref{defPf})
  that 
$({\mathbf P}^{f})^{(-,\epsilon)}
= {\mathbf P}^{(1-\epsilon)f}$. Since 
 $\GM((1-\epsilon)f)=(1-\epsilon)\GM(f)$, 
we see that when $p>0$ and $\pi_{p} \preceq {\mathbf P}^{f},$
$(\pi_{p},{\mathbf P}^{f})$ is downwards movable iff
 it is upwards movable iff $p < \GM(f)$. This implies in particular that
${\mathbf P}^{f}$ is rigid iff $\GM(f)>0$.
} $\Cox$ 
 \end{example}
  The following example is a variant of the one in \cite[Remark 5.4]{LS}.
  \begin{example} \label{ex-epsstab3}
{\rm
  By \cite[Lemma 2.7]{LS}, 
${\mathbf P}^{f} \preceq {\mathbf P}^{g}$ if $f\leq g.$ 
Let $I_{A}$ denote the indicator
  function of some set $A \subseteq [0,1]$ which has 
Lebesgue-measure $1-\delta$ for some $\delta>0.$
  Let $f=1/2 I_{A}$ and $g=3/4 I_{A}.$ There 
exists $\epsilon>0$ such that $f\leq g(1-\epsilon) \leq g,$ 
and so we get that 
${\mathbf P}^{f} \preceq {\mathbf P}^{g(1-\epsilon)} \preceq {\mathbf P}^{g}.$ 
Hence $({\mathbf P}^{f},{\mathbf P}^{g})$
  is downwards movable. However $\GM(g)=0$ which implies 
that ${\mathbf P}^{g}$ does not dominate any product 
  measure with positive density.
  Therefore $({\mathbf P}^{f},{\mathbf P}^{g})$ is not
  upwards movable.
} $\Cox$
  \end{example}

\noindent
For our next example we need the definition of harmonic mean (HM):
\[
  HM(f):=\left(\int_0^1 \frac{dx}{f(x)}  \right)^{-1}.
\]
  \begin{example} \label{ex-epsstab5}
{\rm Let $f(x)=x.$ It is easy to see that $HM(f)=0$ and $GM(f)=1/e>0.$
Since $GM(f)>0$, ${\mathbf P}^{f}$ is rigid. On the other hand,
since $HM(f)=0,$ Theorem 5.16 of \cite{LS} shows that ${\mathbf P}^{f}$ is not uniformly
insertion tolerant.  
} $\Cox$ 
 \end{example}

\section{Uniform insertion tolerance and upwards \\ extractability}
\label{sect:(II)vs(I)}

In this section we focus on uniform
upwards extractability (property (I)) and
uniform insertion tolerance (property (II)). Proposition 
\ref{prop:(II)vs(I)_finite} provides an equivalence between these
properties when $S$ is finite, while Theorem \ref{th:Hajek} exhibits
a contrasting scenario for $S$ countable. 
%,
%and show that for the case of finite $S$,
%(I), (II) and $\supp(\mu)$ is an up-set are equivalent,
%while (II) $\not\Rightarrow$ (I) in the more general countable case,
%including obtaining a translation invariant such example.
\begin{proposition}  \label{prop:(II)vs(I)_finite}
If $S$ is finite and $\mu$ is a probability measure 
on $\{0,1\}^S$, then the following are equivalent:
\begin{description}
\item{\rm (i)} uniform insertion tolerance, 
\item{\rm (ii)} uniform
upwards extractability, and 
\item{\rm (iii)} $\supp(\mu)$ is an up-set. 
\end{description}
\end{proposition}
\begin{theorem}  \label{th:Hajek}
For $S$ countably infinite, there exists a probability measure $\mu$ on 
$\{0,1\}^S$ that is uniformly insertion tolerant but not uniformly upwards extractable. 
Moreover, we can take $\mu$ to be a translation invariant measure on $\{0,1\}^{\mathbb Z}$.
\end{theorem}
{\bf Proof of Proposition \ref{prop:(II)vs(I)_finite}.}
If $\mu$ is uniformly insertion tolerant, then it is immediate that
$\supp(\mu)$ is an up-set. Furthermore uniform upwards 
extractability trivially implies 
uniform insertion tolerance as we have said previously.
We are therefore only left with having to 
show that if $\supp(\mu)$ is an up-set,
then $\mu$ is uniformly upwards extractable.

In what follows, given a configuration $\sigma\in\{0,1\}^S$, $|\sigma|$ will be the number of 1's in $\sigma$.
If there is to exist a $\nu$ such that $\mu=\nu^{(+,\epsilon)}$ with $\epsilon\in[0,1)$, it is not hard 
to see that we must have
\begin{equation} \label{eqnnuformel}
  \nu(\sigma)=\sum_{\tilde{\sigma} \preceq \sigma} (-\epsilon)^{|\sigma|-|\tilde{\sigma}|}
         (1-\epsilon)^{|\tilde{\sigma}|-|S|} \mu(\tilde{\sigma}) \ \forall \sigma \in \{0,1\}^S.
\end{equation}
This can be verified through a direct calculation, but it is easier to calculate
$\nu^{(+,\epsilon)}(\sigma)$ and check that it is indeed 
equal to $\mu(\sigma)$, as follows. 
\begin{eqnarray*}
\lefteqn{\nu^{(+,\epsilon)}(\sigma)} \\
& & = \sum_{\sigma_1 \preceq \sigma} \epsilon^{|\sigma|-|\sigma_1|}(1-\epsilon)^{|S|-|\sigma|}
                                        \nu(\sigma_1) \\
& & =\sum_{\sigma_1 \preceq \sigma} \epsilon^{|\sigma|-|\sigma_1|}(1-\epsilon)^{|S|-|\sigma|}
       \sum_{\sigma_2 \preceq \sigma_1}(-\epsilon)^{|\sigma_1|-|\sigma_2|}
         (1-\epsilon)^{|\sigma_2|-|S|} \mu(\sigma_2) \\
& & =\sum_{\sigma_1 \preceq \sigma} \sum_{\sigma_2 \preceq \sigma_1} \epsilon^{|\sigma|-|\sigma_1|}
          (-\epsilon)^{|\sigma_1|-|\sigma_2|}(1-\epsilon)^{|\sigma_2|-|\sigma|} \mu(\sigma_2) \\
& & = \sum_{\sigma_2}(1-\epsilon)^{|\sigma_2|-|\sigma|} \mu(\sigma_2)
\sum_{\sigma_1:\sigma_2 \preceq \sigma_1 \preceq \sigma}\epsilon^{|\sigma|-|\sigma_1|}(-\epsilon)^{|\sigma_1|-|\sigma_2|}.
\end{eqnarray*}

If we fix $\sigma_2$, then the binomial theorem gives that the last summation is equal to 0
unless $\sigma_2=\sigma$ in which case it is equal to 1. 
We therefore easily obtain that $\nu^{(+,\epsilon)}(\sigma)=\mu(\sigma)$ for every $\sigma.$

What remains is to check that $\nu(\sigma)\geq 0$ for all $\sigma.$ 
From (\ref{eqnnuformel}) it is immediate that $\nu(\sigma)=0$ for 
every $\sigma \not \in \supp(\mu)$ since $\supp(\mu)$ is an up-set. For $\sigma\in\supp(\mu)$ on the 
other hand, it is easy to see that if we do this construction for
different $\epsilon$'s, then we get
\[
\lim_{\epsilon \rightarrow 0} \nu(\sigma) \, = \, \mu(\sigma) \, . 
\]
Since $\mu(\sigma)>0$ for all $\sigma \in \supp(\mu)$ and $|S| <\infty$, for $\epsilon>0$
small enough, we get that $\nu(\sigma)>0$ for all $\sigma \in \supp(\mu).$ This
shows that $\mu$ is $\epsilon$-upwards extractable for all such $\epsilon$. 
$\Cox$

\medskip\noindent
{\bf Proof of Theorem \ref{th:Hajek}.} 
Let $S= \cup_{k=2}^\infty S_k$, where 
\[
S_k= ((k,1), (k,2), \ldots, (k,k)) \, .
\] 
We will take the probability measure
$\mu$ on $\{0,1\}^S$ to be the product measure
\begin{equation}  \label{eq:product_measure}
\mu \, = \, \mu_2 \times \mu_3 \times \cdots
\end{equation}
where each $\mu_k$ is a probability measure on $\{0,1\}^{S_k}$. The 
$\mu_k$'s are constructed as follows, drawing heavily on an example
of Hajek and Berger \cite{HB}. For $\sigma\in\{0,1\}^{S_k}$, let
\begin{equation}  \label{eq:Hajek's_example}
\mu_k(\sigma) = \left\{
\begin{array}{ll}
\frac{4}{3} 2^{-k} & \mbox{if the number of $1$'s in $\sigma$ is even} \\
\frac{2}{3} 2^{-k} & \mbox{if the number of $1$'s in $\sigma$ is odd.} 
\end{array} \right. 
\end{equation}
We may think of $\mu_k$ as the distribution of a $\{0,1\}^{S_k}$-valued random
variable $X_k$ obtained by first tossing a biased coin with
heads-probability $\frac{2}{3}$, and if heads pick the components of
$X_k$ i.i.d.\ $(\frac{1}{2}, \frac{1}{2})$ conditioned on an even number
of $1$'s, while if tails pick the components i.i.d.\ 
$(\frac{1}{2}, \frac{1}{2})$ conditioned on an odd number of $1$'s.
One can also check that this distribution is the same as choosing all but
(an arbitrary) one of the variables according to $\pi_{1/2}$ and then taking the last
variable to be 1 with probability $1/3$ ($2/3$) if there are an even (odd) number of
1's in the other bits. This last description immediately implies that
$\mu_k$ is $\frac{1}{3}$-insertion tolerant. Because
of the product structure in (\ref{eq:product_measure}), this property is
inherited by $\mu$, which therefore is uniformly
insertion tolerant. 

It remains to show that $\mu$ is not uniformly upwards extractable. To this end,
let $X$ be a $\{0,1\}^S$-valued random variable with distribution $\mu$,
and for $k=2,3,\ldots$ let $Y_k$ denote the number of $1$'s
in $X(S_k)$. It is immediate from (\ref{eq:Hajek's_example}) that
\begin{equation}  \label{eq:enhanced_prob_of_even}
\P(Y_k \mbox{ is even}) = {\textstyle \frac{2}{3}}
\end{equation}
for each $k$. Using our last description of $\mu_k$, the weak law of large numbers
implies that
\[
\frac{Y_k}{k} \rightarrow {\textstyle \frac{1}{2}} \mbox{ in probability as $k\to\infty$.}
\]
Hence, in particular,
\begin{equation}  \label{eq:substantial_no_of_0s}
\lim_{k \rightarrow \infty} \P(Y_k \leq k - m) \, = \, 1
\end{equation}
for any fixed $m$. 

Now assume (for contradiction) that $\mu = \nu^{(+, \epsilon)}$ for some
fixed $\epsilon >0$; 
since $\mu$ being $\epsilon_2$-upwards extractable implies it is
$\epsilon_1$-upwards extractable for $\epsilon_1 <\epsilon_2$, we may
without loss of generality assume that $\epsilon \leq 1/3.$
Pick $X'$ according to $\nu$; we may then suppose that
$X$ has been obtained from $X'$ by randomly switching $0$'s to $1$'s 
independently with probability $\epsilon$. The intuition behind the argument leading up 
to a contradiction is that the process of independently flipping $0$'s to $1$'s will cancel 
all preferences of ending up with an even number of $1$'s.

If $X'(S_k)$ contains precisely $l$ $0$'s, then the conditional probability
(given $X'$) that an even number of these switch to
$1$'s when going from $X'$ to $X$ is easily seen to equal
\[
{\textstyle
\frac{1}{2} + \frac{1}{2}(1- 2 \epsilon)^l \, .
}
\]
The easiest way to see this is using an equivalent random mechanism where each 0 
independently ``updates'' with probability $2\epsilon$ and then all the sites 
which have updated then independently actually switch to a 1 with
probability $1/2$.
It follows that the conditional probability (again given $X'$) that $Y_k$ is odd
is at least 
\[
{\textstyle
\min \{ \frac{1}{2} + \frac{1}{2}(1- 2 \epsilon)^l,
\frac{1}{2} - \frac{1}{2}(1- 2 \epsilon)^l \} \, = \, 
\frac{1}{2} - \frac{1}{2}(1- 2 \epsilon)^l. 
}
\]
Now pick $m$ large enough so that 
$\frac{1}{2} - \frac{1}{2}(1- 2 \epsilon)^m > \frac{5}{12}$. Since 
$X' \preceq X$ a.s., we get from (\ref{eq:substantial_no_of_0s}) that
\[
\lim_{k \rightarrow \infty} \P(A_k) \, = \, 1
\]
where $A_k$ is the event that there are at least $m$ $0$'s in
$X'(S_k)$. This gives
\begin{eqnarray*}
\lim_{k \rightarrow \infty} \P(Y_k \mbox{ is odd}) 
& \geq & \lim_{k \rightarrow \infty} \P(Y_k \mbox{ is odd} \, | \, 
A_k) \P(A_k) \\
& \geq & ({\textstyle \frac{1}{2} - \frac{1}{2}(1- 2 \epsilon)^m})
\lim_{k \rightarrow \infty} \P(A_k) \\
& > & {\textstyle \frac{5}{12}} \, . 
\end{eqnarray*}
This clearly contradicts (\ref{eq:enhanced_prob_of_even}).

We now translate this example into the setting of translation 
invariant distributions on $\{0,1\}^{\mathbb Z}$.

Begin with randomly designating either all even integers or all odd integers 
(each with probability $\frac{1}{2}$) in the index set ${\mathbb Z}$
to represent copies of $S_2$. Assume that we happened to choose the
even integers (the other case is handled analogously). Then we
toss another fair coin to decide whether to put i.i.d.\ copies of
$X(S_2)$ on the pairs $\{\ldots, (-4, -2), (0,2), (4,6), \ldots\}$
in ${\mathbb Z}$, 
or on $\{\ldots (-2,0), (2,4), (6,8), \ldots\}$. Then use one more fair coin
to decide whether $\{\ldots -3, 1, 5, 9, \ldots\}$ or
$\{\ldots, -1, 3, 7, 11, \ldots\}$ should be designated for i.i.d.\ copies
of $X(S_3)$, and once this is decided toss a fair three-sided coin to 
choose one of the three possible placements of the length-$3$ blocks in
this subsequence to put these copies. And so on. 

This makes the resulting process $X^*$ translation invariant. Also, since
the property of $\epsilon$-insertion tolerance is obviously closed under
convex combinations, we easily obtain that $X^*$ is $\frac{1}{3}$-insertion 
tolerant
and therefore uniformly insertion tolerant. 

Furthermore, for any $k\geq 2$, 
we may apply (\ref{eq:enhanced_prob_of_even}) to the i.i.d.\ copies
of $X(S_k)$ to deduce that with probability $1$ there will
exist $i \in \{0,1, \ldots, k2^{k-1}-1\}$ such that
\begin{equation}  \label{eq:2/3_limit}
\lim_{n \rightarrow \infty} \frac{1}{n}
\sum_{j=1}^{n} {\bf 1}_{B_{i,j,k}} \, = \, 
{\textstyle \frac{2}{3}}  
\end{equation}
where $B_{i,j,k}$ denotes the event that the number of $1$'s in
\[
\{i+ jk2^{k-1}, i+ jk2^{k-1} + 2^{k-1}, i+ jk2^{k-1} + 2 \cdot 2^{k-1},
\ldots, i+ jk2^{k-1} + (k-1)2^{k-1}\}
\]
is even. The right way to think of $i$ is that it is the first place 
to the right of 
the origin where a copy of $X(S_k)$ starts. The summation variable $j$ 
on the other hand, 
makes us jump to the starting
points of all the other copies of $X(S_k)$ to the right of the origin. 
Furthermore, 
by arguing as in for the non-translation invariant construction, 
we have that if $X^*$ is uniformly upwards extractable, 
then for large $k$ the limit in (\ref{eq:2/3_limit}) will be less than
$1 - \frac{5}{12} = \frac{7}{12}$ for all $i \in \{0,1, \ldots, k2^{k-1}-1\}$.
But this contradicts (\ref{eq:2/3_limit}), so we can conclude that
$X^*$ is not uniformly upwards extractable.
$\Cox$

\medskip\noindent
Note, finally, that the examples in the above proof also show that uniform
finite energy does not imply uniform extractability. 

\section{Rigidity}  \label{sect:slipperiness}

We now proceed to discuss the issue of when a measure is rigid. 
As mentioned in the introduction, any measure
which does not dominate a nontrivial product measure is trivially 
nonrigid and so it would be more interesting to 
have a nonrigid measure which dominates 
a nontrivial product measure; such a measure is provided in
Theorem \ref{extype2_2} below.  

\begin{proposition}  \label{prop:rigid_when_finite}
If $S$ is finite and $\mu$ is a probability measure 
on $\{0,1\}^S$, then the following are equivalent.
\begin{description}
\item{\rm (i)} $\mu$ dominates $\pi_p$ for some $p>0$, 
\item{ \rm (ii)}
$\mu$ is rigid, and 
\item{\rm (iii)} $\mu(\sigma(S)\equiv1)>0$.
\end{description}
\end{proposition}
This does not extend to infinite $S$, as shown in the following result. 

\begin{theorem} \label{extype2_2}
For $S$ countably infinite,
there exists a $\mu$ which dominates a nontrivial product measure $\pi_p$
but is nevertheless nonrigid. 
Moreover, we can take $\mu$ to be a translation invariant measure 
on $\{0,1\}^{\mathbb Z}$.
\end{theorem}
{\bf Proof of Proposition \ref{prop:rigid_when_finite}.}
It is easy to see that the condition that $\mu$ dominates 
$\pi_p$ for some $p>0$
is equivalent to the condition that $\mu(\sigma(S)\equiv 1)>0.$ 
Also, recall that if $\mu$ is rigid it must dominate a 
non-trivial product measure.

To make the proof complete, it only remains to show that (i) and (iii)
of the statement imply that $\mu$ is rigid.
We have $\pi_{p_{\max, \mu}} \preceq \mu$, so that
\begin{equation}  \label{eq:good_ol_dom}
\pi_{p_{\max, \mu}}(A) \, \leq \, \mu(A)
\end{equation}
for all increasing events $A \subseteq \{0,1\}^S$. We next claim that
\begin{equation}  \label{eq:ineq_in_fact_an_equality}
\exists  A\neq\emptyset,\{0,1\}^S \textrm{ such that $A$ is increasing and } \pi_{p_{\max, \mu}}(A)= \mu(A).
\end{equation}
To see this, note that if we had strict inequality
in (\ref{eq:good_ol_dom}) for all such nontrivial increasing events $A$, then
we could find a sufficiently small $\delta>0$ so that
\[
\pi_{p_{\max, \mu}+ \delta}(A) \, < \, \mu(A)
\]
for all such $A$ (this uses the finiteness of $S$), contradicting
the definition of $p_{\max, \mu}$. Now, for such an $A$ we have that
$\mu(A)\geq\mu(\sigma(S)\equiv 1)>0$ and hence for any $\epsilon>0$
\[
\mu^{(-, \epsilon)}(A) \, < \, \mu(A)
\]
(again because $S$ is finite), which in combination with
(\ref{eq:ineq_in_fact_an_equality}) yields
\[
\pi_{p_{\max, \mu}} \, \not\preceq \, \mu^{(-, \epsilon)} \, . 
\]
Since $\epsilon>0$ was arbitrary, $\mu$ is rigid. 
$\Cox$

\medskip\noindent
For the proof of Theorem \ref{extype2_2}, the following elementary
lemma (which is presumably known) is convenient
to have. 
\begin{lemma}  \label{lem:conditioned_binomial}
For $k \geq 1$, $p \in (0,1)$ and $m \in \{0,1, \ldots, k\}$, write
$\rho_{k,p,m}$ for the distribution of a Binomial$(k,p)$ random variable
conditioned on taking value at least $m$. For $p_1 \leq p_2$, we have
\[
\rho_{k,p_1,m} \, \preceq \, \rho_{k,p_2,m} \, . 
\]
\end{lemma}
{\bf Proof.}
For $i=1,2$, let $Y_i$ be a Bin$(k, p_i)$ random variable, and let
$X_i$ be a random variable with distribution $\rho_{k, p_i, m}$. 
What we need to show is that for any $n\in\{m+1, \ldots, k\}$ we have
\[
\frac{\P(X_1 \geq n)}{\P(X_1 < n)} \, \leq \, 
\frac{\P(X_2 \geq n)}{\P(X_2 < n)}
\]
which is the same as showing that
\begin{equation}  \label{eq:need_to_show} 
\frac{\P(X_2 \geq n)}{\P(X_1 \geq n)} \cdot \frac{\P(X_1 < n)}{\P(X_2 < n)}
\, \geq \, 1 \, . 
\end{equation} 
Writing $Z_1$ and $Z_2$ for the probabilities that $Y_1 \geq m$ and
$Y_2 \geq m$, respectively, the left-hand side of (\ref{eq:need_to_show})
becomes
\begin{equation}  \label{eq:first_rewrite}
\frac{\frac{1}{Z_2} \sum_{j=n}^k {k \choose j} 
p_2^j(1-p_2)^{k-j}}{\frac{1}{Z_1} 
\sum_{j=n}^k {k \choose j} p_1^j(1-p_1)^{k-j}} \cdot
\frac{\frac{1}{Z_1} 
\sum_{j=m}^{n-1} {k \choose j} p_1^j(1-p_1)^{k-j}}{\frac{1}{Z_2} 
\sum_{j=m}^{n-1} {k \choose j} 
p_2^j(1-p_2)^{k-j}} \, . 
\end{equation}
Cancelling the $Z_i$'s and introducing the notation 
$\phi_i=\frac{p_i}{1-p_i}$ for $i=1,2$, the expression in
(\ref{eq:first_rewrite}) may further be rewritten as
\begin{eqnarray}  \nonumber
\lefteqn{ \mbox{ } \hspace{-20mm}
\frac{p_2^n(1-p_2)^{k-n} \sum_{j=n}^k {k \choose j} 
\phi_2^{j-n}}{p_1^n(1-p_1)^{k-n} \sum_{j=n}^k {k \choose j} \phi_1^{j-n}}
\cdot
\frac{p_1^n(1-p_1)^{k-n} \sum_{j=m}^{n-1} {k \choose j} 
\phi_1^{j-n}}{p_2^n(1-p_2)^{k-n} \sum_{j=m}^{n-1} {k \choose j} 
\phi_2^{j-n}} = } \\
& = &
\frac{ \sum_{j=n}^k {k \choose j} 
\phi_2^{j-n}}{ \sum_{j=n}^k {k \choose j} \phi_1^{j-n}}
\cdot
\frac{ \sum_{j=m}^{n-1} {k \choose j} 
\phi_1^{j-n}}{ \sum_{j=m}^{n-1} {k \choose j} 
\phi_2^{j-n}} \, .  
\label{eq:second_rewrite}
\end{eqnarray}
Note now that $\phi_1 \leq \phi_2$, so that
\[
\sum_{j=n}^k {k \choose j} 
\phi_2^{j-n} \, \geq \, \sum_{j=n}^k {k \choose j} \phi_1^{j-n}
\]
and
\[
\sum_{j=m}^{n-1} {k \choose j} 
\phi_1^{j-n} \, \geq \, \sum_{j=m}^{n-1} {k \choose j} 
\phi_2^{j-n} \, . 
\]
Hence, the expression in (\ref{eq:second_rewrite}) is greater than or equal
to $1$, so (\ref{eq:need_to_show}) is verified and the lemma is
established. $\Cox$

\medskip\noindent
{\bf Proof of Theorem \ref{extype2_2}.}
As in the proof of Theorem \ref{th:Hajek}, we take 
$S= \cup_{k=2}^\infty S_k$ where 
$S_k= ((k,1), (k,2), \ldots, (k,k))$, and 
the probability measure
$\mu$ on $\{0,1\}^S$ to be the product measure
\[
\mu \, = \, \mu_2 \times \mu_3 \times \cdots
\]
where each $\mu_k$ is a probability measure on $\{0,1\}^{S_k}$. This 
time, we take the $\mu_k$'s to be as follows. For $\sigma \in \{0,1\}^{S_k}$,
set
\begin{equation}  \label{eq_def_of_blocks}
\mu_k(\sigma) = \left\{
\begin{array}{ll}
k^{-1}2^{-k} & \mbox{if the number of $1$'s in $\sigma$ is exactly $1$}, \\
1 - 2^{-k} & \mbox{if } \sigma=(1,1,1,\ldots, 1),  \\
0 & \mbox{otherwise}. 
\end{array}  \right. 
\end{equation}
We now make three claims about the $\mu_k$ measures:
\begin{description}
\item{\sc Claim 1.} $p_{\max, \mu_k} \geq \frac{1}{2}$ for all $k$.
\item{\sc Claim 2.} $\lim_{k \rightarrow\infty}p_{\max, \mu_k}= \frac{1}{2}$.
\item{\sc Claim 3.} For any fixed
$\epsilon < \frac{1}{2}$, we have for all $k$ sufficiently large that
\[
\mu_k^{(-,\epsilon)} \succeq \pi_{\frac{1}{2}}
\]
where $\pi_{\frac{1}{2}}$ is product measure with $p=\frac{1}{2}$ on
$\{0,1\}^{S_k}$.
\end{description}
We slightly postpone proving the claims, and first show how they
imply the existence of a nonrigid measure that dominates $\pi_{\frac{1}{2}}$. 

Let us modify $S$ and $\mu$ slightly by setting, for $m \geq 2$, 
\[
\tilde{S}_m \, = \, \cup_{k=m}^\infty S_k
\]
and
\begin{equation}  \label{eq:projected_product_structure}
\tilde{\mu}_m \, = \, \mu_m \times \mu_{m+1} \times \cdots
\end{equation}
so that in other words $\tilde{\mu}_m$ is the probability measure on
$\{0,1\}^{\tilde{S}_m}$ which arises by projecting $\mu$ on 
$\{0,1\}^{\tilde{S}_m}$. 

Using the product structure (\ref{eq:projected_product_structure}), we
get from {\sc Claim 1} that $p_{\max, \tilde{\mu}_m} \geq \frac{1}{2}$
(for any $m$), and from {\sc Claim 2} that 
$p_{\max, \tilde{\mu}_m} \leq \frac{1}{2}$ (for any $m$). Hence
\[
p_{\max, \tilde{\mu}_m} \, = \, {\textstyle \frac{1}{2} }
\]
for any $m$. Fixing $\epsilon \in (0,1/2)$, we can also deduce from 
(\ref{eq:projected_product_structure}) and {\sc Claim 3} that
\begin{equation}  \label{eq:punchline_giving_nonrigidity}
\tilde{\mu}_m^{(-, \epsilon)} \, \succeq \, \pi_{\frac{1}{2}}
\, = \, \pi_{p_{\max, \tilde{\mu}_m}}
\end{equation}
for $m$ sufficiently large. For such $m$ we thus have that $\tilde{\mu}_m$
is nonrigid. 

It remains to prove {\sc Claim 1}, {\sc Claim 2} and {\sc Claim 3}. 

{\sc Claim 1} is the same as saying that $\mu_k \succeq 
\pi_{\frac{1}{2}}$. This
is immediate to verify, but the best way to think about it is as follows. 
Suppose that we pick $X_k \in \{0,1\}^{S_k}$ according to 
$\pi_{\frac{1}{2}}$,
and if $X_k= (0,0, \ldots, 0)$ then we switch one of the $0$'s (chosen
uniformly at random) to a $1$, while otherwise we switch {\em all} $0$'s
to $1$'s. The resulting random element of $\{0,1\}^{S_k}$ then
has distribution $\mu_k$. 

To prove {\sc Claim 2}, it suffices (in view of {\sc Claim 1}) to prove that
\[
\limsup_{k \rightarrow\infty}p_{\max, \mu_k}\le \frac{1}{2}
\]
and to this end it is enough to show for any $\delta>0$ that
\begin{equation} \label{eq:doesnt_dom_1/2+delta}
\mu_k \, \not\succeq \pi_{\frac{1}{2}+\delta}
\end{equation}
for all sufficiently large $k$. Let $A_k$ denote the event of seeing
at most one $1$ in $\{0,1\}^{S_k}$; then $A_k$ is a decreasing event and
its complement $\neg A_k$ is increasing. Now simply note that
\begin{equation}  \label{eq:first_ratio_of_probs}
\frac{\mu_k(A_k)}{\pi_{\frac{1}{2}+\delta}(A_k)} \, = \, 
\frac{(\frac{1}{2})^k}{(\frac{1}{2}- \delta)^k + 
k (\frac{1}{2}+ \delta)(\frac{1}{2}- \delta)^{k-1}}
\end{equation}
which tends to $\infty$ as $k \rightarrow \infty$. Hence, taking $k$ 
large enough gives $\mu_k(A_k) > \pi_{\frac{1}{2}+\delta}(A_k)$, so that
$\mu_k(\neg A_k) < \pi_{\frac{1}{2}+\delta}(\neg A_k)$ and
(\ref{eq:doesnt_dom_1/2+delta}) is established, proving {\sc Claim 2}. 

To prove {\sc Claim 3}, note first that both $\pi_{\frac{1}{2}}$ and
$\mu_k^{(-, \epsilon)}$ are invariant under permutations of $S_k$, 
so that it suffices to show for $k$ large that
\begin{equation}  \label{eq:consider_only_Bn}
\mu_k^{(-, \epsilon)} (B_n) \, \leq \, \pi_{\frac{1}{2}} (B_n)
\end{equation}
for $n=0,1, \ldots, k-1$, where $B_n$ is the event of seeing at most
$n$ $1$'s in $S_k$. For $n=0$ we get
\begin{equation}  \label{eq:second_ratio_of_probs}
\frac{\mu_k^{(-, \epsilon)} (B_0)}{\pi_{\frac{1}{2}} (B_0)} \, = \, 
\frac{(\frac{1}{2})^k \epsilon + (1- (\frac{1}{2})^k ) 
\epsilon^k}{(\frac{1}{2})^k}
\end{equation}
while for $n=1$
\begin{equation}  \label{eq:third_ratio_of_probs}
\frac{\mu_k^{(-, \epsilon)} (B_1)}{\pi_{\frac{1}{2}} (B_1)} \, = \,
\frac{ (\frac{1}{2})^k + (1- (\frac{1}{2})^k ) 
(\epsilon^k + k \epsilon^{k-1}(1-\epsilon))}{(k+1)(\frac{1}{2})^k } \, .
\end{equation}
The right-hand sides of (\ref{eq:second_ratio_of_probs}) and 
(\ref{eq:third_ratio_of_probs})
tend to $\epsilon$ and $0$, respectively, 
as $k \rightarrow \infty$, so (\ref{eq:consider_only_Bn})
is verified for $n=0$ and $1$ (and $k$ large enough). 
To verify (\ref{eq:consider_only_Bn}) for $n\geq 2$ (and all such $k$), 
define two random variables $Y$ and $Y'$ as the number of $1$'s in
two random elements of $\{0,1\}^{S_k}$ with respective distributions
$\mu_k^{(-, \epsilon)}$ and $\pi_{\frac{1}{2}}$. Note that $Y$ conditioned
on taking value at least $2$ has the same distribution as a 
Bin $(k, 1- \epsilon)$ random variable conditional on taking value
at least $2$, while the conditional distribution of $Y'$ given that it
is at least $2$, is that of a Bin $(k, \frac{1}{2})$ variable conditioned
on being at least $2$. Defining $\rho_{k, (1-\epsilon), 2}$ and
$\rho_{k, \frac{1}{2}, 2}$ as in Lemma \ref{lem:conditioned_binomial},
we thus have for $n \in \{2, \ldots, k-1\}$ that
\begin{equation}  \label{eq:first_compl_prob}
\mu_k^{(-, \epsilon)} (B_n) \, = \, 1 - (1-\mu_k^{(-, \epsilon)} (B_1))
(1- \rho_{k, (1-\epsilon), 2}(B_n))
\end{equation}
and
\begin{equation}  \label{eq:second_compl_prob}
\pi_{\frac{1}{2}} (B_n) \, = \, 1 - (1-\pi_{\frac{1}{2}} (B_1))
(1- \rho_{k, \frac{1}{2}, 2}(B_n)) \, . 
\end{equation}
But we have already seen that 
$\mu_k^{(-, \epsilon)} (B_1) \leq \pi_{\frac{1}{2}} (B_1)$, and Lemma
\ref{lem:conditioned_binomial} tells us that 
$\rho_{k, (1-\epsilon), 2}(B_n) \leq \rho_{k, \frac{1}{2}, 2}(B_n)$, so
(\ref{eq:first_compl_prob}) and (\ref{eq:second_compl_prob}) yield
\[
\mu_k^{(-, \epsilon)} (B_n) \, \leq \, \pi_{\frac{1}{2}} (B_n) \, ,
\]
and {\sc Claim 3} is established. 

Finally, we translate this example into the setting of translation 
invariant distributions on $\{0,1\}^{\mathbb Z}$.
The measure $\tilde{\mu}_m$ 
can be turned into a translation invariant measure $\tilde{\mu}^*_m$
on $\{0,1\}^{\mathbb Z}$ 
by the same independent-copy-and-paste procedure as in
Theorem \ref{th:Hajek}. The property
\[
\pi_{\frac{1}{2}} \, \preceq \, (\tilde{\mu}^*_m)^{(-, \epsilon)}
\]
is obviously inherited from (\ref{eq:punchline_giving_nonrigidity}). Thus,
in order to show that $\tilde{\mu}^*_m$ is nonrigid, it only remains to
show that it does not stochastically dominate $\pi_{\frac{1}{2}+ \delta}$
for any $\delta >0$. This follows using (\ref{eq:first_ratio_of_probs})
by an argument analogous to (\ref{eq:2/3_limit}) in
Theorem \ref{th:Hajek}: If we pick $k$ depending on $\delta$ as in
the justification of {\sc Claim 2}, then, under $\tilde{\mu}^*_m$,  
certain infinite arithmetic progressions will have subsequences 
of length $k$ which contain at most one $1$ often
enough (under spatial averaging) that the corresponding event
has $\pi_{\frac{1}{2}+ \delta}$-measure $0$. We omit the details. 
$\Cox$

\begin{remark}
\label{rem:neither(II)or(III)}
{\rm
The measure $\tilde{\mu}_m$ is obviously not uniformly insertion 
tolerant, and we have thus demonstrated the existence of a measure for which 
property (IV) holds while neither (II) nor (III) does. 
}
$\Cox$
\end{remark}
\begin{remark}
{\rm For any $p \in (0,1)$, the construction above can be modified
by replacing $2^{-k}$ by $p^k$ in (\ref{eq_def_of_blocks}). Proceeding
as in the rest of the proof yields the result that for any
$p, \epsilon\in (0,1)$ such that $p+\epsilon<1$, there exists a measure
$\mu$ on $\{0,1\}^S$ where $S$ is countably infinite, with
the property that
$p_{\max, \mu} = p$ and 
\[
\pi_{p_{\max, \mu}} \preceq \mu^{(-, \epsilon)} \, . 
\]
This is obviously sharp.
$\Cox$
}
\end{remark}

  \section{Further results on rigidity} \label{secFKG}

In this section, we continue the study of rigidity, and prove 
Theorem \ref{thmFKG} and Proposition \ref{prop:neg}.

\medskip\noindent
The proof of  Theorem \ref{thmFKG} will make use of the following
technical lemma.
  \begin{lemma} \label{lemmaFKG}
  Let $\mu$ be a measure on $\{0,1\}^{{\mathbb Z}^{d}}.$ Assume that it is 
  $\delta$-insertion tolerant for some $\delta >0$. If for some $p\in (0,1)$ and $\epsilon >0$
  \begin{equation} \label{eqnassump}
    \mu^{(-,\epsilon)}(\sigma(\{1,\ldots,n\}^d)\equiv 0)\leq (1-p)^{n^d} \textrm{ for all } n\geq 0,
  \end{equation}
  then there exists $p'>p$ such that
  \[
    \mu(\sigma(\{1,\ldots,n\}^d)\equiv 0)\leq (1-p')^{n^d} \textrm{ for all } n\geq 0.
  \]
  \end{lemma}
  {\bf Proof.}
  Let $X\sim \mu$ and $Z\sim \pi_{1-\epsilon}$ be independent and let $X^{(-,\epsilon)}=
  \min(X,Z).$ It is easy to see using the $\delta$-insertion tolerance that for any $s \in \{1,\ldots,n\}^d,$ 
  and any $\zeta \in \{0,1\}^{\{1,\ldots,n\}^d \setminus s}$
  \begin{eqnarray*}
    \lefteqn{{\mathbb P}(X(s)=1 \cap X(\{1,\ldots,n\}^d \setminus s)\equiv \zeta)} \\
    & &   \geq \frac{\delta}{1-\delta}{\mathbb P}(X(s)=0 \cap X(\{1,\ldots,n\}^d \setminus s)\equiv \zeta).
  \end{eqnarray*}
  Iterating this, we get that for any $\xi \in \{0,1\}^{\{1,\ldots,n\}^d}$
  \[
    {\mathbb P}(X(\{1,\ldots,n\}^d)\equiv \xi) \geq \left ( \frac{\delta}{1-\delta} \right ) ^{|\xi|} 
                {\mathbb P}(X(\{1,\ldots,n\}^d)\equiv 0).
  \]  
  Here $|\xi|$ denotes the cardinality of the set $\{s\in \{1,\ldots,n\}^d:\xi(s)=1\}.$
  We get that
\begin{eqnarray*}
\lefteqn{{\mathbb P}(X^{(-,\epsilon)}(\{1,\ldots,n\}^d)\equiv 0)} \\
&=& \sum_{\xi \in \{0,1\}^{\{1,\ldots,n\}^d}}{\mathbb P}
(X^{(-,\epsilon)}(\{1,\ldots,n\}^d)\equiv 0|X(\{1,\ldots,n\}^d)\equiv \xi)\\ 
& & \times {\mathbb P}(X(\{1,\ldots,n\}^d)\equiv \xi) \\
&\geq & \sum_{\xi \in \{0,1\}^{\{1,\ldots,n\}^d}}{\mathbb P}
(X^{(-,\epsilon)}(\{1,\ldots,n\}^d)\equiv 0|X(\{1,\ldots,n\}^d)\equiv \xi) \\
& & \times \left ( \frac{\delta}{1-\delta} \right ) ^{|\xi|}{\mathbb P}
(X(\{1,\ldots,n\}^d)\equiv 0) \\
&=&  \sum_{\xi \in \{0,1\}^{\{1,\ldots,n\}^d}}\epsilon^{|\xi|} 
\left ( \frac{\delta}{1-\delta} \right ) ^{|\xi|}{\mathbb P}
(X(\{1,\ldots,n\}^d)\equiv 0) \\
&=&   \left (1+ \frac{\epsilon \delta}{1-\delta} \right )^{n^d}{\mathbb P}
(X(\{1,\ldots,n\}^d)\equiv 0).
  \end{eqnarray*}
  Therefore if (\ref{eqnassump}) holds we can conclude that
  \[
    \mu(\sigma(\{1,\ldots,n\}^d)\equiv 0)\leq \left (\frac{1-\delta}{1-\delta+\epsilon \delta} \right)^{n^d}
                                      (1-p)^{n^d} \, , 
  \]
and we are done.
$\Cox$

\medskip  \noindent
  {\bf Proof of Theorem \ref{thmFKG}.}
  The case $p_{\max,\mu}=1$ is trivial and we therefore assume that $p_{\max,\mu}\in(0,1).$
  In \cite{LigS}, it is shown that if $\mu$ is downwards FKG and if 
  \begin{equation} \label{liggett}
   \mu(\sigma(\{1,\ldots,n\}^d)\equiv 0) \leq (1-p)^{n^d} \textrm{ for all } n\geq 0,
  \end{equation}
 then $\pi_p \preceq \mu.$ Therefore if $\pi_{p_{\max,\mu}} \preceq \mu^{(-,\epsilon)}$ for some
$\epsilon >0$, then (\ref{eqnassump}) trivially holds (with $p=p_{\max,\mu}$)
and so we can conclude from Lemma~\ref{lemmaFKG} and the above result in \cite{LigS}
that $\pi_{p'} \preceq \mu$ for some $p'> p_{\max,\mu}$, a contradiction.
$\Cox$

\medskip  \noindent
We now define conditional negative association.

  \begin{definition}
  A probability measure $\mu$ on $\{0,1\}^{\mathbb Z}$ is said to have 
{\bf conditional negative association} if for any finite $S \subset {\mathbb Z}$ 
and any two increasing functions $f,g$ that are measurable on disjoint
subsets of ${\mathbb Z}\setminus S$,
  \[
    \mu(fg|\sigma(S)) \leq   \mu(f|\sigma(S))\mu(g|\sigma(S)).
  \] 
  \end{definition}
  We will use the fact (see \cite{LS}) that for conditionally 
negatively associated measures $\mu$, we have
$\pi_{\rho} \preceq \mu$ iff 
  \begin{equation} \label{eqncNA}
    \mu(\sigma(\{1,\ldots,n\})\equiv 1)\geq {\rho}^n \ \forall \ n\geq 1. 
  \end{equation}

\noindent
{\bf Proof of Proposition \ref{prop:neg}.} 
We note that the case $p_{\max,\mu}=1$ is trivial and we 
  therefore assume that $p_{\max,\mu}\in(0,1).$  Assume that
  $\pi_{p_{\max,\mu}} \preceq \mu^{(-,\epsilon)}$ for some $\epsilon > 0$.
We then get that 
  \[
p_{\max,\mu}^n \leq \mu^{(-,\epsilon)}(\sigma(\{1,\ldots,n\})\equiv 1)=(1-\epsilon)^n\mu(\sigma(\{1,\ldots,n\})\equiv 1).
  \]
  Hence 
  \[
    \mu(\sigma(\{1,\ldots,n\})\equiv 1) \geq   
\left( \frac{p_{\max,\mu}}{1-\epsilon} \right)^n.
  \]
  Therefore $\pi_{\frac{p_{\max,\mu}}{1-\epsilon}} \preceq \mu$ by the 
result quoted in connection with (\ref{eqncNA}). 
This is a contradiction since $p_{\max,\mu} < \frac{p_{\max,\mu}}{1-\epsilon}$.
$\Cox$

  \section{Proof of main result} \label{types}

  \begin{lemma} \label{lemmapos}
If $\mu$ is uniform upwards extractable, then
for any $\epsilon>0$ there exists a 
  $\delta>0$ such that $(\mu^{(-,\epsilon)})^{(+,\delta)} \preceq \mu.$
  \end{lemma}
  {\bf Proof.} Let $\nu$ and $\alpha>0$ be such that $\mu=\nu^{(+,\alpha)}.$ 
One can easily compute that for any $\alpha$,
$\epsilon$, and $\delta$, we have that
  \[
    ((\mu^{(+,\alpha)})^{(-,\epsilon)})^{(+,\delta)}
           =\mu^{(-,\epsilon(1-\delta),+,\alpha(1-\epsilon)+\alpha \epsilon \delta+(1-\alpha)\delta)}.
  \]
  Now, given $\epsilon>0,$ choose $\delta>0$ such that 
  $\alpha(1-\epsilon)+\alpha \epsilon \delta+(1-\alpha)\delta<\alpha.$ 
We therefore get that 
  \[
    (\mu^{(-,\epsilon)})^{(+,\delta)}=((\nu^{(+,\alpha)})^{(-,\epsilon)})^{(+,\delta)} 
     \preceq \nu^{(-,\epsilon(1-\delta),+,\alpha)} \preceq \nu^{(+,\alpha)}=\mu.
  \]
$\Cox$

  \begin{lemma} \label{lemmatype1}
  Given a probability measure $\mu$ on $\{0,1\}^{S},$ assume that for 
every $\epsilon>0$, there exists a $\delta>0$ such that 
$(\mu^{(-,\epsilon)})^{(+,\delta)} \preceq \mu$. Then $\mu$ is rigid.
  \end{lemma}
  {\bf Proof.}
  The case $p_{\max,\mu}=1$ is trivial, and we will therefore assume 
that  $p_{\max,\mu}\in (0,1)$. Assume for contradiction that $\mu$ is 
nonrigid. Then there exists an $\epsilon>0$ such that 
  $\pi_{p_{\max,\mu}} \preceq \mu^{(-,\epsilon)}$. By assumption there 
exists a $\delta>0$ such that 
 $(\mu^{(-,\epsilon)})^{(+,\delta)} \preceq \mu$. 
Hence $(\pi_{p_{\max,\mu}})^{(+,\delta)} \preceq 
  (\mu^{(-,\epsilon)})^{(+,\delta)} \preceq \mu$. 
Since $p_{\max,\mu}<1,$ $(\pi_{p_{\max,\mu}})^{(+,\delta)}$ is
  a product measure with density strictly larger than $p_{\max,\mu}.$ 
  This is a contradiction.
$\Cox$

\medskip\noindent
We remark that we do not know whether the reverse statement of 
Lemma \ref{lemmatype1} is true. It would also be interesting to know if
the sufficient condition in this lemma follows from uniform insertion
tolerance.
 
Example \ref{ex-epsstab5} provides us with an example of a $\mu$ which is 
on one hand rigid but on the other hand not uniformly
insertion tolerant. However, since it relies heavily on results not presented
in this paper, we give here another more ``hands on'' example.
It is a variant of \cite[Remark 6.4]{LS} and 
shows that the reverse statement of Lemma \ref{lemmapos} is false.

  \begin{example} \label{ex-epsstab4}
{\rm
  Let $\{X_{i}\}_{i \in {\mathbb N}}$ be defined in the following way. 
For every even $i \geq 0$, let independently $(X_{i},X_{i+1})$ be $(1,1)$ 
or $(0,0)$ with probability $1/2$ each. Let $\mu_e$ denote the
  distribution of this process. For $\epsilon, \delta>0$ let 
$\{X_{i}^{(-,\epsilon(1-\delta),+,\delta)}\}_{i \in {\mathbb N}}$
 be a sequence of random variables with distribution 
$\mu_e^{(-,\epsilon(1-\delta),+,\delta)}=(\mu_e^{(-,\epsilon)})^{(+,\delta)}$.
By noting that for any $\epsilon>0$ there exists a $\delta>0$ such that for even $i$
  \[
    {\mathbb P}(\max(X_{i}^{(-,\epsilon(1-\delta),+,\delta)},
X_{i+1}^{(-,\epsilon(1-\delta),+,\delta)})=1)<{\textstyle \frac{1}{2}},
  \]
  we see that for the same choice of $\epsilon,\delta$ we get that 
$(\mu_e^{(-,\epsilon)})^{(+,\delta)} \preceq \mu_e$. Lemma \ref{lemmatype1}
  gives us that $\mu_e$ is rigid.
  However, it is easy to see that $\mu_e$ is not uniform insertion tolerant.

  The only drawback with this construction is that it is not translation 
invariant. However this is easily fixed. Let $\mu_o$ be the distribution 
of $\{X_{i+1}\}_{i \in {\mathbb N}},$ i.e. it is $\mu_e$ shifted over by 1.
  Define the measure $\mu$ by 
  \[
 {\textstyle   \mu=\frac{1}{2}\mu_e+\frac{1}{2}\mu_o.  }
  \]
 It is easy to check that
  \[
 {\textstyle   (\mu^{(-,\epsilon)})^{(+,\delta)}=\frac{1}{2}(\mu_e^{(-,\epsilon)})^{(+,\delta)}+\frac{1}{2}(\mu_o^{(-,\epsilon)})^{(+,\delta)}
          \preceq \frac{1}{2}\mu_e+\frac{1}{2}\mu_o=\mu.  }
  \]
By Lemma \ref{lemmatype1}, it follows that $\mu$ is rigid.
On the other hand, clearly
$$
\mu(\sigma(0)=1|\sigma(1)=0,\sigma(2)=\sigma(3)=1)=0
$$
and hence $\mu$ is not uniformly insertion tolerant.
} $\Cox$
\end{example}

\noindent
  {\bf Proof of Theorem \ref{thm1}.}
  Lemma \ref{lemmapos} together with Lemma \ref{lemmatype1} 
shows that  property (I) implies property (III) and all the other
implications were indicated in the introduction. As far as all of the
reversed implications claimed not to hold, we continue as follows.
Example \ref{ex-epsstab4} together with Lemma \ref{lemmatype1}
(or example \ref{ex-epsstab5})
shows that (III) does not imply (II) (and hence that
(III) does not imply (I) and that (IV) does not imply (II)).
Theorem \ref{extype2_2} implies that (IV) does not imply (III). 
Finally, Theorem \ref{th:Hajek} shows that (II) does not imply (I). 
Also, all of these examples were translation invariant measures on
$\{0,1\}^{{\mathbb Z}}.$
$\Cox$

\medskip
\noindent
{\bf Proof of Theorem \ref{thm:finite}.}
This follows immediately from Propositions \ref{prop:(II)vs(I)_finite}
and \ref{prop:rigid_when_finite}.
$\Cox$

\medskip
\noindent
We feel, finally, that it is worth mentioning the following result,
which is an easy consequence of Lemma \ref{lemmapos}.

  \begin{corollary}
  Assume that $(\mu_1,\mu_2)$ is downwards movable and that $\mu_2$
is uniformly upwards extractable. Then
$(\mu_1,\mu_2)$ is also upwards movable.
  \end{corollary}

\end{document}